\documentclass{article}
\usepackage{amssymb}
\usepackage{amsmath,amsthm}
\usepackage{amsfonts}
\usepackage{graphicx}
\usepackage{graphics}

\vspace{.9cm}

\newtheorem{Theorem}{Theorem}[section]

\newtheorem{cor}{Corollary}[section]
\newtheorem{lemma}{Lemma}[section]

\begin{document}

\date{\small\textsl{\today}}
\title{ Strong convergence of a new composite iterative method for equilibrium problems and fixed point problems in Hilbert spaces}
\author{  \large Vahid Darvish
          and ~\large S. M. Vaezpour
\footnote{Corresponding author.\newline {\em  E-mail addresses:}
\newline vahiddarvish@aut.ac.ir (V. Darvish),
\newline vaez@aut.ac.ir (S. M. Vaezpour).}
           $\vspace{.2cm} $  \\
\small{\em Department of Mathematics and Computer Science,}\vspace{-1mm}\\
\small{\em Amirkabir University of Technology,}\vspace{-1mm}\\
\small{\em Hafez Ave., P.O. Box 15875-4413, Tehran, Iran }\\}
\maketitle

\begin{abstract} In this paper, first we introduce a new mapping
for finding a common fixed point of an infinite family of
nonexpansive mappings then we consider iterative method for finding
a common element of the set of fixed points of an infinite family of
nonexpansive mappings, the set of solutions of an equilibrium
problem and the set of solutions of the variational inequality for
$\alpha$-inverse-strongly monotone mapping in a Hilbert space. We
show that under suitable conditions, the sequence converges strongly
to a common element of the above three sets. Our results presented
in this paper improve and extend other results.

\medskip


\medskip


\end{abstract}

\section{Introduction and preliminaries}
Throughout this paper, we always assume that $H$ is a real Hilbert
space with inner product, $\langle .,. \rangle$ and norm, $\|.\|$,
respectively and $C$ is a nonempty closed convex subset of $H$. We use
F(T) to denote the set of fixed point of $T$. Recall that a mapping
$T:C\longrightarrow C$ is called nonexpansive if
$\|Tx-Ty\|\leq\|x-y\|, \ \forall x,y\in C$, and is a contraction, if
there exists a constant $\alpha\in(0,1)$ such that
$\|Tx-Ty\|\leq\alpha\|x-y\|, \ \forall x,y \in C$.\\
Let $B:C\longrightarrow H$ be a mapping, the classical variational
inequality problem is to find a $u \in C$ such that,
\begin{equation}
\langle Bu,v-u\rangle \geq0, \ \forall v\in C.\label{classical}
\end{equation}
The set of solution of variational inequality (\ref{classical}) is
denoted by VI(B,C).\\
An operator $A$ is said to be strongly positive if there exists a
constant $\bar{\gamma}>0$ such that, $\langle
Ax,x\rangle\geq\bar{\gamma}\|x\|^{2}, \ \forall x\in H$. It is easy
to see that the following is true, $$u\in VI(B,C)\Longleftrightarrow
u=P_{C}(u-\lambda Bu),  \ \lambda > 0.$$ Let $F$ be a bifunction of
$C\times C \longrightarrow \mathbb{R}$. The equilibrium problem for
$F:C\times C\longrightarrow \mathbb{R}$ is to find $x\in C$ such
that,
\begin{eqnarray}
F(x,y)\geq 0 \label{equilibrium}
\end{eqnarray}
for all $y\in C$, the set of solution of (\ref{equilibrium}) is
denoted by EP(F). It is well known that for every point $x\in H$,
there exists a unique nearest point in $C$, denoted by $P_{C}x$,
such that $\|x-P_{C}x\|\leq\|x-y\|, \ \forall y\in C$.\\
$P_{C}$ is called the metric projection of $H$ onto $C$ and is a
nonexpansive mapping which satisfies $$\langle x-y,P_{C}x-P_{C}y
\rangle\geq \|P_{C}x-P_{C}y\|^{2}, \ \forall x,y\in H.$$ Moreover,
$P_{C}x$ is characterized by the following properties, $P_{C}x\in C$
and $$\langle x-P_{C}x,y-P_{C}y\rangle\leq0,$$
$$\|x-y\|^{2}\geq\|x-P_{C}x\|^{2}+\|y-P_{C}x\|^{2},$$ for all $x\in
H$, $y\in C$.\\

\noindent For solving the equilibrium problem for a bifunctional $F:C\times C
\longrightarrow \mathbb{R}$, let us assume that $F$ satisfies the
following conditions:\\
(A1) $F(x,x)=0$, for all $x\in C$;\\
(A2) $F$ is monotone, i.e., $F(x,y)+F(y,x)\leq 0$, for all $x,y \in
C$; \\
(A3) for each $x,y,z \in C$, $\lim_{t\to 0}F(tz+(1-t)x,y)\leq
F(x,y)$;\\
(A4) for each $x\in C$, $y\longmapsto F(x,y)$ is convex and lower
semicontinuous.\\

\noindent Recall that a mapping $B:C\longrightarrow H$ is said to be:\\
(i) monotone if $\langle Bu-Bv,u-v \rangle \geq0, \ \forall u,v\in
C$;\\
(ii) $L$-Lipschitz if there exists a constant $L>0$ such that
$\|Bu-Bv\|\leq L\|u-v\|,$ for all $u,v \in C$;\\
(iii) $\alpha$-inverse-strongly monotone, if there exists a positive
real number $\alpha$ such that \begin{equation} \langle
Bu-Bv,u-v\rangle\geq
\alpha\|Bu-Bv\|^{2}, \ \forall u,v \in C \label{alphainverse}.\\
\end{equation}
We know that any $\alpha$-inverse-strongly monotone mapping $B$ is
monotone and $(\frac{1}{\alpha}$)-Lipschitz continuous.\\

\noindent Numerous problems in physics, optimization and economics
are to find a solution of (\ref{equilibrium}). Let $A$ be strongly
positive bounded linear operator on $H$, a typical problem is to
minimize a quadratic function over the set of the fixed points of a
nonexpansive mapping on a real Hilbert space
H, i.e.,\\
$$\min _{x\in F(T)}\frac{1}{2}\langle Ax,x\rangle -\langle x,b\rangle
,$$ where $b$ is a given point in $H$.\\

\noindent Takahashi \cite{tak} introduced an iterative method for
finding a common element of equilibrium points of $F$ and the set of
fixed points of a nonexpansive mapping in Hilbert spaces as follows:

\begin{eqnarray*}
&&F(u_{n},y)+\frac{1}{r_{n}}\langle y-u_{n},u_{n}-x_{n}\rangle\geq0,
\
\forall y\in C,\\
&&x_{n+1}=\alpha_{n}f(x_{n})+(1-\alpha_{n})Tu_{n}.\\
\end{eqnarray*}
They proved under certain conditions $\{x_{n}\}$ and $\{u_{n}\}$
converges strongly to $z\in F(T)\bigcap EP(F)$, where $z \in
P_{F(T)\bigcap EP(F)}f(z).$\\

\noindent Moreover, Plubtieng and Pungpaeng \cite{plub} introduced
the following iterative scheme:
\begin{eqnarray*}
&&F(u_{n},y)+\frac{1}{r_{n}}\langle y-u_{n},u_{n}-x_{n}\rangle
\geq0,
\ \forall y\in H,\\
&&x_{n+1}=\alpha_{n}\gamma f(x_{n})+(I-\alpha_{n}A)Tu_{n}, \ \forall
n\in \mathbb{N}.\\
\end{eqnarray*}
They proved under certain appropriate conditions on $\{\alpha_{n}\}$
and $\{r_{n}\}$ the sequence $\{x_{n}\}$ and $\{u_{n}\}$ generated
by above scheme converges strongly to the unique solution of the
variational inequality $\langle (A-\gamma f)z,x-z \rangle \geq0, \
\forall x\in F(T)\bigcap EP(F).$ \\

\noindent Zhao and He \cite{zhao} introduced the following process:\\
\begin{eqnarray*}
&&F(u_{n},y)+\frac{1}{r_{n}}\langle y-u_{n},u_{n}-x_{n}\rangle
\geq0, \
\forall y\in C,\\
&&y_{n}=s_{n}P_{C}(u_{n}-\lambda_{n}Au_{n})+(1-s_{n})x_{n},\\
&&x_{n+1}=\alpha_{n}u+\beta_{n}x_{n}+\gamma_{n}W_{n}(P_{C}(y_{n}-\lambda_{n}Ay_{n})).\\
\end{eqnarray*}
They proved under some conditions the sequences $\{x_{n}\}$ and
$\{u_{n}\}$ converges strongly to\\
$z\in\bigcap_{i=1}^{\infty}F(T_{i})\bigcap VI(A,C)\bigcap EP(F),$
where $W_{n}$ is defined by Shimoji and Takahashi \cite{shi}, as
follows:
\begin{eqnarray}
U_{n,n+1}&=&I,\nonumber\\
U_{n,n}&=&\gamma_{n}T_{n}U_{n,n+1}+(1-\gamma_{n})I,\nonumber\\
U_{n,n-1}&=&\gamma_{n-1}T_{n-1}U_{n,n}+(1-\gamma_{n-1})I,\nonumber\\
&\vdots&\nonumber\\
U_{n,k}&=&\gamma_{k}T_{k}U_{n,k+1}+(1-\gamma_{k})I,\nonumber\\
U_{n,k-1}&=&\gamma_{k-1}T_{k-1}U_{n,k}+(1-\gamma_{k-1})I,\nonumber\\
&\vdots&\nonumber\\
U_{n,2}&=&\gamma_{2}T_{2}U_{n,3}+(1-\gamma_{2})I,\nonumber\\
W_{n}&=&U_{n,1}=\gamma_{1}T_{1}U_{n,2}+(1-\gamma_{1})I,\nonumber\label{12}\\
\end{eqnarray}
where $\gamma_{1},\gamma_{2},...$ are real numbers such that $0\leq
\gamma_{n}\leq1, T_{1},T_{2},...$ are an infinite family of mappings
of $H$ into itself. Note that the nonexpansivity of each $T_{i}$
ensures the nonexpansivity of $W_{n}$. \\

\noindent Qin, Cho and Kang \cite{qin} introduced the following
iterative
scheme:\\
\begin{eqnarray*}
x_{1}&=&x\in H, \ \ chosen \ \ arbitrary,\\
z_{n}&=&\lambda_{n}x_{n}+(1-\lambda_{n})W_{n}x_{n},\\
y_{n}&=&\beta_{n}\gamma f(z_{n})+(I-\beta_{n}A)z_{n},\\
x_{n+1}&=&\alpha_{n}x_{n}+(1-\alpha_{n})y_{n}, \ \forall n\geq1.\\
\end{eqnarray*}
\begin{equation}\label{corq}
\end{equation}
They proved that under certain assumptions on the sequences
$\{\alpha_{n}\}$, $\{\beta_{n}\}$ and $\{\lambda_{n}\}$, $\{x_{n}\}$
converges strongly to a common fixed point of the infinite family of
nonexpansive mappings.\\

Throughout this paper, inspired by Qin, Cho and Kang \cite{qin} and
Zhao and He \cite{zhao}, we introduce a composite iterative method
for
infinite family of nonexpansive mappings as follows:\\
\begin{eqnarray}
&&F(u_{n},y)+\frac{1}{r_{n}}\langle y-u_{n},u_{n}-x_{n}\rangle\geq0,
\
\forall y\in C,\nonumber \\
&&z_{n}=\lambda_{n}u_{n}+(1-\lambda_{n})W_{n}u_{n}, \nonumber\\
&&y_{n}=\beta_{n}\gamma
f(z_{n})+(I-\beta_{n}A)P_{C}(z_{n}-\gamma_{n}Bz_{n}), \nonumber \\
&&x_{n+1}=\alpha_{n}x_{n}+(1-\alpha_{n})y_{n}, \nonumber \label{main} \\
\end{eqnarray}
where $W_{n}$ is defined as (\ref{12}), $f$ is a contraction on $H$,
$A$ is strongly positive bounded linear self-adjoint operator, $B$
is $\alpha$-inverse-strongly monotone. We prove under certain
appropriate assumptions on the sequences $\{\alpha_{n}\}$,
$\{\beta_{n}\}$, $\{\lambda_{n}\}$ and $\{\gamma_{n}\}$, the
sequence $\{x_{n}\}$ defined by (\ref{main}), converges strongly to
a common element of the set of fixed points of a finite family of
nonexpansive mappings, the set of solutions of an equilibrium
problem and the set of solutions of the variational inequality for
$\alpha$-inverse-strongly monotone mapping.\\
We need the following lemmas for the proof of our main results.

\begin{lemma}\label{r} The following inequality holds in a Hilbert space
H,
$$\|x+y\|^2 \leq \|x\|^2 + 2\langle y,x+y\rangle, \forall x,y \in H.$$
\end{lemma}

\begin{lemma}\label{1} \cite{xu} Assume $\{\alpha_n\}$ is a sequence of nonnegative real
numbers such that $\alpha_{n+1} \leq (1-\gamma_n)\alpha_n +
\delta_n, \ \ n\geq 1$ , where $\{\gamma_n \}$ is a sequence in
$(0,1)$ and $\delta_n$ is a
sequence in $\mathbb{R}$ such that:\\
\begin{enumerate}
\item $\sum_{n=1}^{\infty}\gamma_{n}=\infty ,$
\item $\limsup_{n \to \infty}(\frac{\delta_{n}}{\gamma_{n}})\leq 0
$ or $\sum_{n=1}^{\infty}|\delta_{n}|< \infty , $
\end{enumerate}
then $\lim_{n\to\infty}\alpha_{n}=0.$
\end{lemma}
\begin{lemma}\label{d} \cite{marino} Assume that A is a strongly positive linear bounded
self-adjoint operator on a Hilbert space $H$ with coefficient
$\bar{\gamma}$ and $0<\rho\leq \|A\|^{-1}$, then $\|I-\rho
A\|\leq1-\rho\bar{\gamma}$.
\end{lemma}
\begin{lemma}\label{55}\cite{suzuki}
Let $\{x_{n}\}$ and $\{y_{n}\}$ be bounded sequences in Banach space
$X$ and let $\{\beta_{n}\}$ be a sequence in $[0,1]$ with
$0<\liminf_{n\to\infty}\beta_{n}\leq\limsup_{n\to\infty}\beta_{n}<1$.
Suppose that $x_{n+1}=(1-\beta_{n})y_{n}+\beta_{n}x_{n}$ for all
integers $n\geq0$ and
$$\limsup_{n\to\infty}(\|y_{n+1}-y_{n}\|-\|x_{n+1}-x_{n}\|)\leq0.$$
Then, $\lim_{n\to\infty}\|y_{n}-x_{n}\|=0$.\\
\end{lemma}
\begin{lemma}\label{66}\cite{blum}
Let $C$ be a nonempty closed convex subset of $H$ and let $F$ be a
bifunction of $C\times C$ into $R$ satisfying (A1)-(A4). Let $r>0$
and $x\in H$. Then, there exists $z\in C$ such that,
$$F(z,y)+\frac{1}{r}\langle y-z,z-x\rangle\geq0, \ \forall y\in C.$$
\end{lemma}
\begin{lemma}\label{77}\cite{comb}
Assume that $F\colon C\times C\longrightarrow R$ satisfies
(A1)-(A4). For $r>0$ and $x\in H$, define a mapping
$T_{r}:H\longrightarrow C$ as
follows:\\
$T_{r}(x)=\{z\in C : F(z,y)+\frac{1}{r}\langle y-z,z-x\rangle\geq0,
\ \forall y\in C\}$ for all $z\in H$. Then the following hold:\\
\begin{enumerate}
\item $T_{r}$ is single-valued;
\item $T_{r}$ is firmly nonexpansive, i.e., for any $x, y\in H$,
$\|T_{r}x-T_{r}y\|^{2}\leq\langle T_{r}x-T_{r}y,x-y\rangle$;
\item $F(T_{r})=EP(F)$;
\item $EP(F)$ is closed and convex.
\end{enumerate}
\end{lemma}
\begin{lemma}\cite{shi} Let $C$ be nonempty closed convex subset of a Hilbert
space, let $T_{i}:C\longrightarrow C$ be an infinite family of
nonexpansive mappings with $\bigcap_{i=1}^{\infty}F(T_{i})\neq
\emptyset$ and let ${\gamma_{i}}$ be a real sequence such that
$0<\gamma_{i}\leq\gamma<1$ for all $i\geq1$ then ,
\begin{enumerate}
\item $W_{n}$ is nonexpansive and $F(W_{n})=\bigcap_{i=1}^{n}F(T_{i})$ for
each $n\geq 1$ ,
\item For each $x\in C$ and for each positive integer $k$, the
$\lim_{n\to\infty}U_{n,k}$ exists ,
\item The mapping $W:C\longrightarrow C$ defined by :
$$Wx:=\lim_{n\to\infty}W_{n}x=\lim_{n\to\infty}U_{n,1}x \  \ \ x\in
C,$$ \\

is a nonexpansive mapping satisfying
$F(W)=\bigcap_{i=1}^{\infty}F(T_{i})$ and is called the $W$-mapping
generated by $T_{1},T_{2},...$ and $\gamma_{1},\gamma_{2},... \ .$
\end{enumerate}
\end{lemma}
\begin{lemma}\cite{shi}\label{2568}
Let $C$ be a nonempty closed convex subset of a Hilbert space $H$,
let $T_{i}:C\longrightarrow C$ be an infinite family of nonexpansive
mappings with $\bigcap_{i=1}^{\infty}F(T_{i})\neq\emptyset$ and let
$\gamma_{i}$ be a real sequence such that $0<\gamma_{i}\leq\gamma<1$
for all $i\geq1$, if $K$ is any bounded subset of $C$ then,
$$\limsup_{n\to\infty}\|Wx-W_{n}x\|=0, \ \ \ x\in K.$$
\end{lemma}

\section{Main Results}
In this section, we prove strong convergence theorem.\\
\begin{Theorem}\label{maintheorem} Let $H$ be a real Hilbert space. Let $F$ be a
bifunction from $C\times C \longrightarrow \mathbb{R}$ satisfying
(A1)-(A4) and let $B$ be an $\alpha$-inverse-strongly monotone
mapping of $C$ into $H$, and let $\{T_{i}:C\longrightarrow C\}$ be a
infinite family of nonexpansive mappings with
$F:=\bigcap_{i=1}^{\infty}F(T_{i})\bigcap VI(B,C)\bigcap EP(F)\neq
\emptyset$. Suppose $A$ is a strongly positive linear bounded
self-adjoint operator with the coefficient $\bar{\gamma}\geq0$ and
$0<\gamma\leq\frac{\bar{\gamma}}{\alpha}$. Let $\{\alpha_{n}\}$,
$\{\beta_{n}\}$, $\{\lambda_{n}\}$ and $\{\delta_{n}\}$ be sequences
in $[0,1]$ satisfying the following conditions:\\
(C1) $\sum_{n=0}^{\infty}\beta_{n}=\infty$,
$\lim_{n\to\infty}\beta_{n}=0$;\\
(C2) $\sum_{i=1}^{\infty}|\lambda_{n}-\lambda_{n+1}|<\infty$;\\
(C3)
$\sum_{i=1}^{\infty}|\alpha_{n}-\alpha_{n+1}|<\infty;$\\
(C4) $\exists \lambda\in[0,1]; \lambda_{n}<\lambda, \ \forall
n\geq1;$\\
(C5) $\sum_{n=1}^{\infty}|\gamma_{n}-\gamma_{n+1}|<\infty,
\gamma_{n}\in[a,b], a,b\in(0,2\alpha);$\\
(C6) $\liminf_{n\to\infty}r_{n}>0$ and $\sum_{n=1}^{\infty}|r_{n+1}-r_{n}|<\infty$.\\
Then the sequences $\{x_{n}\}$ defined by (\ref{main}) converges
strongly to $q\in F$, where $q=P_{F}(\gamma f+(I-A))(q)$ which
solves the following variational inequality: $$\langle \gamma
f(q)-Aq,p-q\rangle\leq0, \ \forall p\in F.$$\\
\end{Theorem}
\noindent{proof}\ : For all $x,y\ \in C$ and $\gamma_{n}\in
(0,2\alpha)$, we note that

\begin{eqnarray}
\|(I-\gamma_{n}B)x-(I-\gamma_{n}B)y\|^{2}&=&\|(x-y)-\gamma_{n}(Bx-By)\|^{2}\nonumber\\
&=&\|x-y\|^{2}-2\gamma_{n}\langle x-y,Bx-By\rangle
+\gamma_{n}^{2}\|Bx-By\|^{2}\nonumber\\
&\leq&\|x-y\|^{2}+\gamma_{n}(\gamma_{n}-2\alpha)\|Bx-By\|^{2}\nonumber\\
&\leq& \|x-y\|^{2}\nonumber\label{111},\\
\end{eqnarray}
which implies that $I-\gamma_{n}B$ is nonexpansive. Noticing that
$A$ is a linear bounded self-adjoint operator so we have,
$$\|A\|=sup\{|\langle Ax,x\rangle|: x\in H, \|x\|=1\}.$$
Since $\beta_{n}\to 0$ as ${n\to\infty}$, by (C1) we may assume
without loss of generality $\beta_{n}<\|A\|^{-1}, \ \forall
n\in\mathbb{N}$.\\

$$\langle(I-\beta_{n}Ax),x\rangle=1-\beta_{n}\langle Ax,x\rangle
\geq1-\beta_{n}\|A\|\geq0,$$\\
So $I-\beta_{n}A$ is positive.\\
It follows that,\\
\begin{eqnarray*}
\|I-\beta_{n}A\|&=&\sup\{\langle(I-\beta_{n}A)x,x\rangle : x\in H,
\|x\|=1\}\\
&=&\sup\{1-\beta_{n}\langle Ax,x\rangle : x\in H, \|x\|=1\}\\
&\leq& 1-\beta_{n}\bar{\gamma}.\\
\end{eqnarray*}
Let $p\in\bigcap_{i=1}^{\infty}F(T_{i})\bigcap VI(B,C)\bigcap EP(F)$
and $\{T_{r_{n}}\}$ be a sequence of mapping defined as in lemma
(\ref{77}), then $p=P_{C}(p-\gamma_{n}Bp)=T_{r_{n}}p$. Put
$t_{n}=P_{C}(z_{n}-\gamma_{n}Bz_{n})$, from (\ref{111}) we have,\\
\begin{eqnarray*}
\|t_{n}-p\|&=&\|P_{C}(z_{n}-\gamma_{n}Bz_{n})-P_{C}(p-\gamma_{n}Bp)\|\\
&\leq&\|z_{n}-\gamma_{n}Bz_{n}-p+\gamma_{n}Bp\|\\
&\leq&\|z_{n}-p\|.\\
\end{eqnarray*}
Next we show that $\{x_{n}\}$ is bounded. Let $p\in F$, from the
definition of $T_{r}$, we know that $u_{n}=T_{r_{n}}x_{n}$, it
follows that
$\|u_{n}-p\|=\|T_{r_{n}}x_{n}-T_{r_{n}}p\|\leq\|x_{n}-p\|$.
Furthermore we have,\\
\begin{eqnarray*}
\|z_{n}-p\|&=&\|\lambda_{n}u_{n}+(1-\lambda_{n})W_{n}u_{n}-p\|\\
&=&\|\lambda_{n}(u_{n}-p)+(1-\lambda_{n})(W_{n}u_{n}-p)\|\\
&\leq&\lambda_{n}\|u_{n}-p\|+(1-\lambda_{n})\|W_{n}u_{n}-p\|\\
&=&\|u_{n}-p\|\leq\|x_{n}-p\|,\\
\end{eqnarray*}
and hence,\\
\begin{eqnarray*}
\|y_{n}-p\|&=&\|\beta_{n}\gamma
f(z_{n})+((I-\beta_{n}A))P_{C}(z_{n}-\gamma_{n}Bz_{n})-p\|\\
&=&\|\beta_{n}(\gamma
f(z_{n})-Ap)+(I-\beta_{n}A)(P_{C}(z_{n}-\gamma_{n}Bz_{n})-p)\|\\
&\leq&\beta_{n}\|\gamma
f(z_{n})-Ap\|+(1-\beta_{n}\bar{\gamma})\|z_{n}-p\|\\
&\leq&\beta_{n}\gamma\|f(z_{n})-f(p)\|+\beta_{n}\|\gamma
f(p)-Ap\|+(1-\beta_{n}\bar{\gamma})\|z_{n}-p\|\\
&\leq&\beta_{n}\gamma\alpha\|z_{n}-p\|+\beta_{n}\|\gamma
f(p)-Ap\|+(1-\beta_{n}\bar{\gamma})\|z_{n}-p\|\\
&\leq&\beta_{n}\gamma\alpha\|x_{n}-p\|+\beta_{n}\|\gamma
f(p)-Ap\|+(1-\beta_{n}\bar{\gamma})\|x_{n}-p\|\\
&=&\left(1-\beta_{n}(\bar{\gamma}-\gamma\alpha)\right)\|x_{n}-p\|+\beta_{n}\|\gamma
f(p)-Ap\|.\\
\end{eqnarray*}
It follows that,\\
\begin{eqnarray*}
\|x_{n+1}-p\|&=&\|\alpha_{n}x_{n}+(1-\alpha_{n})y_{n}-p\|\\
&\leq&\alpha_{n}\|x_{n}-p\|+(1-\alpha_{n})\|y_{n}-p\|\\
&\leq&\alpha_{n}\|x_{n}-p\|+(1-\alpha_{n})[(1-\beta_{n}(\bar{\gamma}-\gamma\alpha))\|x_{n}-p\|+\beta_{n}\|\gamma
f(p)-Ap\|]\\
&=&[1-\beta_{n}(\bar{\gamma}-\gamma\alpha)(1-\alpha_{n})]\|x_{n}-p\|+\beta_{n}(\bar{\gamma}-\gamma\alpha)(1-\alpha_{n})\frac{\|\gamma
f(p)-Ap\|}{\bar{\gamma}-\gamma\alpha},\\
\end{eqnarray*}
by simple induction we have,\\
$$\|x_{n}-p\|\leq\max\{\|x_{1}-p\|,\frac{\|Ap-\gamma
f(p)\|}{\bar{\gamma}-\gamma\alpha}\}.$$ Which gives that the
sequence
$\{x_{n}\}$ is bounded, so $\{y_{n}\}$ and $\{z_{n}\}$.\\

\noindent Next, we claim that
$\lim_{n\to\infty}\|x_{n+1}-x_{n}\|=0$.\\
First we have $u_{n}=T_{r_{n}}x_{n}$, $u_{n-1}=T_{r_{n}}x_{n-1}$,
so\\
\begin{equation}\label{el}
F(u_{n-1},y)+\frac{1}{r_{n-1}}\langle
y-u_{n-1},u_{n-1}-x_{n-1}\rangle \geq0, \ \forall y\in C,
\end{equation}
and \\
\begin{equation}\label{jav}
F(u_{n},y)+\frac{1}{r_{n}}\langle y-u_{n},u_{n}-x_{n}\rangle\geq0, \
\forall y\in C,\\
\end{equation}
putting $y=u_{n}$ in (\ref{el}) and $y=u_{n-1}$ in (\ref{jav}), we
have,\\
$$F(u_{n-1},u_{n})+\frac{1}{r_{n-1}}\langle
u_{n}-u_{n-1},u_{n-1}-x_{n-1}\rangle \geq0,$$
and\\
$$F(u_{n},u_{n-1})+\frac{1}{r_{n}}\langle
u_{n-1}-u_{n},u_{n}-x_{n}\rangle\geq0.$$ Hence, $$\langle
u_{n}-u_{n-1},\frac{u_{n-1}-x_{n-1}}{r_{n-1}}-\frac{u_{n}-x_{n}}{r_{n}}\rangle\geq0,$$
therefore we have,\\
$$\langle u_{n}-u_{n-1},u_{n-1}-u_{n}+u_{n}-x_{n-1}-\frac{r_{n-1}}{r_{n}}(u_{n}-x_{n})\rangle\geq0.$$
We can assume that there is a real number $b$ such that $r_{n}>b>0$,
for all $n\in \mathbb{N}$. Then we have,\\
\begin{eqnarray*}
\|u_{n}-u_{n-1}\|^{2}&\leq&\langle u_{n}-u_{n-1},
x_{n}-x_{n-1}+(1-\frac{r_{n-1}}{r_{n}})(u_{n}-x_{n})\rangle\\
&\leq&
\|u_{n}-u_{n-1}\|[\|x_{n}-x_{n-1}\|+|1-\frac{r_{n-1}}{r_{n}}|\|u_{n}-x_{n}\|].\\
\end{eqnarray*}
So,
\begin{eqnarray*}
\|u_{n}-u_{n-1}\|&\leq&\|x_{n}-x_{n-1}\|+\frac{1}{r_{n}}|r_{n}-r_{n-1}|\|u_{n}-x_{n}\|\\
&\leq&\|x_{n}-x_{n-1}\|+\frac{1}{b}|r_{n}-r_{n-1}|L,\\
\end{eqnarray*}
where $L=\sup\{\|u_{n}-x_{n}\|: n\in \mathbb{N}\}$.\\
On the other hand,\\
\begin{eqnarray*}
\|z_{n}-z_{n-1}\|&=&\|\lambda_{n}u_{n}+(1-\lambda_{n})W_{n}u_{n}-\lambda_{n-1}u_{n-1}-(1-\lambda_{n-1})W_{n-1}u_{n-1}\|\\
&\leq&\lambda_{n}\|u_{n}-u_{n-1}\|+|\lambda_{n}-\lambda_{n-1}|\|u_{n-1}-W_{n-1}u_{n-1}\|\\
&&+(1-\lambda_{n})\|W_{n}u_{n}-W_{n}u_{n-1}\|+(1-\lambda_{n})\|W_{n}u_{n-1}-W_{n-1}u_{n-1}\|\\
&\leq&\|u_{n}-u_{n-1}\|+|\lambda_{n}-\lambda_{n-1}|\|u_{n-1}-W_{n-1}u_{n-1}\|+(1-\lambda_{n})\|W_{n}u_{n-1}-W_{n-1}u_{n-1}\|.\\
\end{eqnarray*}
Also by nonexpansivity of $T_{i}$ and the definition of $W_{n}$ we
have,\\
\begin{eqnarray*}
\|W_{n}u_{n-1}-W_{n-1}u_{n-1}\|&=&\|\mu_{1}T_{1}U_{n,2}u_{n-1}-(1-\mu_{1})u_{n-1}-\mu_{1}T_{1}U_{n-1,2}u_{n-1}-(1-\mu_{1})u_{n-1}\|\\
&\leq&\mu_{1}\|U_{n,2}u_{n-1}-U_{n-1,2}u_{n-1}\|\\
&=&\mu_{1}\|\mu_{2}T_{2}U_{n,3}u_{n-1}-(1-\mu_{2})u_{n-1}-\mu_{2}T_{2}U_{n-1,3}u_{n-1}-(1-\mu_{2})u_{n-1}\|\\
&\leq&\mu_{1}\mu_{2}\|U_{n,3}u_{n-1}-U_{n-1,3}u_{n-1}\|\\
&\vdots&\\
&\leq&(\prod_{i=1}^{n-1}\mu_{i})\|U_{n,n}u_{n-1}-U_{n-1,n}u_{n-1}\|\\
&\leq& M_{1}(\prod_{i=1}^{n-1}\mu_{i}),\\
\end{eqnarray*}
where $M_{1}\geq0$ is a constant such that
$\|U_{n,n}u_{n-1}-U_{n-1,n}u_{n-1}\|\leq M_{1}.$\\
Now we have,\\
$\|z_{n}-z_{n-1}\|=\|u_{n}-u_{n-1}\|+|\lambda_{n}-\lambda_{n-1}|\|u_{n-1}-W_{n-1}u_{n-1}\|+(1-\lambda_{n})M_{1}.$\\
Also,\\
\begin{eqnarray*}
\|y_{n}-y_{n-1}\|&=&\|\beta_{n}\gamma
f(z_{n})+(I-\beta_{n}A)t_{n}-\beta_{n-1}\gamma
f(z_{n-1})-(I-\beta_{n-1}A)t_{n-1}\|\\
&\leq&\beta_{n}\gamma\alpha\|z_{n}-z_{n-1}\|+\gamma|\beta_{n}-\beta_{n-1}|\|f(z_{n-1})\|+(1-\beta_{n}\bar{\gamma})\|t_{n}-t_{n-1}\|+|\beta_{n-1}-\beta_{n}|\|At_{n-1}\|,\\
\end{eqnarray*}
where,\\
\begin{eqnarray*}
\|t_{n}-t_{n-1}\|&=&\|P_{C}(z_{n}-\gamma_{n}Bz_{n})-P_{C}(z_{n-1}-\gamma_{n-1}Bz_{n-1})\|\\
&\leq&\|(z_{n}-z_{n-1})-(\gamma_{n}Bz_{n}-\gamma_{n-1}Bz_{n-1})\|\\
&\leq&\|z_{n}-z_{n-1}\|+\frac{\gamma_{n}}{\alpha}\|z_{n}-z_{n-1}\|+|\gamma_{n}-\gamma_{n-1}|\|Bz_{n-1}\|.\\
\end{eqnarray*}
So,\\
\begin{eqnarray*}
\|y_{n}-y_{n-1}\|&\leq&\beta_{n}\gamma\alpha\|z_{n}-z_{n-1}\|+\gamma|\beta_{n}-\beta_{n-1}|\|f(z_{n-1})\|+|\beta_{n-1}-\beta_{n}|\|At_{n-1}\|\\
&&+(1-\beta_{n}\bar{\gamma})\left((\frac{\alpha+\gamma_{n}}{\alpha})\|z_{n}-z_{n-1}\|+|\gamma_{n}-\gamma_{n-1}|\|Bz_{n-1}\|\right)\\
\end{eqnarray*}
\begin{eqnarray*}
&\leq&\left(\beta_{n}\gamma\alpha+(1-\beta_{n}\bar{\gamma})(\frac{\alpha+\gamma_{n}}{\alpha})\right)\|z_{n}-z_{n-1}\|\\
&&+\left(|\gamma_{n}-\gamma_{n-1}|(1-\beta_{n}\bar{\gamma})+|\beta_{n}-\beta_{n-1}|\right)M_{2},\\
\end{eqnarray*}
where, $M_{2}=\sup\{\gamma\|f(z_{n-1}\|+\|Bz_{n-1}\|+\|At_{n-1}\|\}$
\begin{eqnarray*}
&\leq&(1-\beta_{n}(\bar{\gamma}-\gamma\alpha)+\frac{\gamma_{n}}{\alpha})[\|u_{n}-u_{n-1}\|\\
&&+|\lambda_{n}-\lambda_{n-1}|\|u_{n-1}-W_{n-1}u_{n-1}\|+(1-\lambda_{n})M_{1}\prod_{i=1}^{n-1}\mu_{i}]\\
&&+p_{n}M_{2},\\
\end{eqnarray*}
where,
$p_{n}=(|\gamma_{n}-\gamma_{n-1}|(1-\beta_{n}\bar{\gamma})+|\beta_{n}-\beta_{n-1}|).$\\
\begin{eqnarray*}
&\leq&(1-\beta_{n}(\bar{\gamma}\alpha-\gamma\alpha)+\frac{\gamma_{n}}{\alpha})[\|x_{n}-x_{n-1}\|+\frac{1}{b}|r_{n}-r_{n-1}|L]\\
&&+(1+\frac{\gamma_{n}}{\alpha})[|\lambda_{n}-\lambda_{n-1}|\|u_{n-1}-W_{n-1}u_{n-1}\|\\
&&+(1-\lambda_{n})M_{1}\prod_{i=1}^{n-1}\mu_{i}]+p_{n}M_{2}.\\
\end{eqnarray*}
We note that,\\
\begin{eqnarray*}
\|x_{n+1}-x_{n}\|&=&\|\alpha_{n}x_{n}+(1-\alpha_{n})y_{n}-\alpha_{n-1}x_{n-1}-(1-\alpha_{n-1})y_{n-1}\|\\
&=&\|\alpha_{n}x_{n}-\alpha_{n}x_{n-1}+\alpha_{n}x_{n-1}-\alpha_{n-1}x_{n-1}+(1-\alpha_{n})y_{n}-(1-\alpha_{n})y_{n-1}\\
&&+(1-\alpha_{n})y_{n-1}-(1-\alpha_{n-1})y_{n-1}\|\\
&\leq&\alpha_{n}\|x_{n}-x_{n-1}\|+|\alpha_{n}-\alpha_{n-1}|\|x_{n-1}\|+(1-\alpha_{n})\|y_{n}-y_{n-1}\|+|\alpha_{n-1}-\alpha_{n}|\|y_{n-1}\|.\\
\end{eqnarray*}
If we put $K=\sup\{\|x_{n-1}\|+\|y_{n-1}\|\}.$\\
So,\\
$\|x_{n+1}-x_{n}\|\leq\alpha_{n}\|x_{n}-x_{n-1}\|+2K|\alpha_{n}-\alpha_{n-1}|+(1-\alpha_{n})\|y_{n}-y_{n-1}\|.$
By substitute $\|y_{n}-y_{n-1}\|$ in above inequality we have,\\
\begin{eqnarray*}
\|x_{n+1}-x_{n}\|&\leq&\alpha_{n}\|x_{n}-x_{n-1}\|+(1-\alpha_{n})(1-\beta_{n}(\bar{\gamma}-\gamma\alpha)+\frac{\gamma_{n}}{\alpha})[\|x_{n}-x_{n-1}\|\\
&&+\frac{1}{b}|r_{n}-r_{n-1}|L]\\
&&+(1-\alpha_{n})(1+\frac{\gamma_{n}}{\alpha})[\|u_{n-1}-W_{n-1}u_{n-1}\||\lambda_{n}-\lambda_{n-1}|\\
&&+(1-\lambda_{n})M_{1}]+(1-\alpha_{n})p_{n}M_{2}+2K|\alpha_{n}-\alpha_{n-1}|\\
&=&(3-\beta_{n}(\bar{\gamma}-\gamma\alpha))\|x_{n}-x_{n-1}\|+s_{n},\\
\end{eqnarray*}
where,\\
\begin{eqnarray*}
s_{n}&=&\{(1-\alpha_{n})(1-\beta_{n}(\bar{\gamma}-\gamma\alpha)+\frac{\gamma_{n}}{\alpha})(\frac{1}{b}|r_{n}-r_{n-1}|L)\\
&&+(1-\alpha_{n})(1+\frac{\gamma_{n}}{\alpha})[\|u_{n-1}-W_{n-1}u_{n-1}\||\lambda_{n}-\lambda_{n-1}|+(1-\lambda_{n})M_{1}\prod_{i=1}^{n-1}\mu_{i}]\\
&&+(1-\alpha_{n})p_{n}M_{2}+2K|\alpha_{n}-\alpha_{n-1}|\}.\\
\end{eqnarray*}
So,\\
$\|x_{n+1}-x_{n}\|\leq(1-\beta_{n}(\frac{\bar{\gamma}-\gamma\alpha}{3}))\|x_{n}-x_{n-1}\|+s_{n},$
 and by lemma (\ref{1}), we have $$\|x_{n}-x_{n-1}\|\to0.$$
\begin{equation}\label{234}
\end{equation}
Now we show that $\lim_{n\to\infty}\|y_{n}-u_{n}\|=0$. We have,\\
\begin{eqnarray*}
\|y_{n+1}-y_{n}\|&=&\|(I-\beta_{n+1}A)(t_{n+1}-t_{n})-(\beta_{n+1}-\beta_{n})At_{n}\\
&&+\gamma[\beta_{n+1}(f(z_{n+1})-f(z_{n}))+f(z_{n})(\beta_{n+1}-\beta_{n})]\|\\
&\leq&(1-\beta_{n+1}\bar{\gamma})\|t_{n+1}-t_{n}\|+|\beta_{n+1}-\beta_{n}|\|At_{n}\|\\
&&+\gamma\beta_{n+1}\alpha\|z_{n+1}-z_{n}\|+|\gamma_{n+1}-\gamma_{n}|\|Bz_{n}\|\\
&&+\gamma|\beta_{n+1}-\beta_{n}|\|f(z_{n})\|\\
&\leq&(1-\beta_{n+1}(\bar{\gamma}-\gamma\alpha))\|z_{n+1}-z_{n}\|+|\gamma_{n+1}-\gamma_{n}|\|Bz_{n}\|+|\beta_{n+1}-\beta_{n}|M_{3}\\
&\leq&\|z_{n+1}-z_{n}\|+|\gamma_{n+1}-\gamma_{n}|\|Bz_{n}\|+|\beta_{n+1}-\beta_{n}|M_{3},\\
\end{eqnarray*}
where $M_{3}=\sup\{\|At_{n}\|+\gamma\|f(z_{n})\|\}$. \\
So,\\
\begin{eqnarray*}
\|y_{n+1}-y_{n}\|&\leq&\|u_{n+1}-u_{n}\|+\|u_{n}-W_{n}u_{n}\||\lambda_{n+1}-\lambda_{n}|\\
&&+(1-\lambda_{n+1})M_{1}\prod_{i=1}^{n-1}\mu_{i}+|\gamma_{n+1}-\gamma_{n}|\|Bz_{n}\|+|\beta_{n+1}-\beta_{n}|M_{3},\\
\end{eqnarray*}
which implies,\\
\begin{eqnarray*}
\|y_{n+1}-y_{n}\|-\|u_{n+1}-u_{n}\|&\leq&\|u_{n}-W_{n}u_{n}\||\lambda_{n+1}-\lambda_{n}|+(1-\lambda_{n+1})M_{1}\prod_{i=1}^{n-1}\mu_{i}\\
&&+|\beta_{n+1}-\beta_{n}|M_{3}+|\gamma_{n+1}-\gamma_{n}|\|Bz_{n}\|.\\
\end{eqnarray*}
By conditions (C1), (C2) and (C5) we have,
$\limsup\{\|y_{n+1}-y_{n}\|-\|u_{n+1}-u_{n}\|\}\leq0$  and by lemma
(\ref{55}),
we have $\limsup_{n\to\infty}\|y_{n}-u_{n}\|=0$.\\
\begin{equation}\label{333}
\end{equation}
Also we show that $\lim_{n\to\infty}\|x_{n}-u_{n}\|=0$. We have,\\
\begin{eqnarray*}
\|x_{n+1}-p\|^{2}&=&\|\alpha_{n}x_{n}+(1-\alpha_{n})y_{n}-p\|^{2}\\
&\leq&\alpha_{n}\|x_{n}-p\|^{2}+(1-\alpha_{n})\|y_{n}-p\|^{2}\\
&\leq&\alpha_{n}\|x_{n}-p\|^{2}+(1-\alpha_{n})(\|y_{n}-u_{n}\|+\|u_{n}-p\|)^{2}\\
&=&\alpha_{n}\|x_{n}-p\|^{2}+(1-\alpha_{n})\|y_{n}-u_{n}\|^{2}+2(1-\alpha_{n})\|y_{n}-u_{n}\|\|u_{n}-p\|+(1-\alpha_{n})\|u_{n}-p\|^{2}.\\
\end{eqnarray*}
Also,\\
\begin{eqnarray}
\|u_{n}-p\|^{2}&=&\|T_{r_{n}}x_{n}-T_{r_{n}}p\|^{2}\leq\langle
T_{r_{n}}x_{n}-T_{r_{n}}p,x_{n}-p\rangle\nonumber\\
&=&\langle u_{n}-p,x_{n}-p\rangle =
\frac{1}{2}(\|u_{n}-p\|^{2}+\|x_{n}-p\|^{2}-\|x_{n}-u_{n}\|^{2}).\nonumber\label{222}\\
\end{eqnarray}
Now by (\ref{222}) we have,\\
\begin{eqnarray*}
\|x_{n+1}-p\|^{2}&\leq&\alpha_{n}\|x_{n}-p\|^{2}+(1-\alpha_{n})(\|x_{n}-p\|^{2}-\|x_{n}-u_{n}\|^{2})+d_{n}\\
&=&\|x_{n}-p\|^{2}-(1-\alpha_{n})\|x_{n}-u_{n}\|^{2}+d_{n},\\
\end{eqnarray*}
where,
$$d_{n}=(1-\alpha_{n})\|y_{n}-u_{n}\|^{2}+2(1-\alpha_{n})\|y_{n}-u_{n}\|\|u_{n}-p\|.$$
Therefore,\\
\begin{eqnarray*}
(1-\alpha_{n})\|x_{n}-u_{n}\|^{2}&\leq&\|x_{n}-p\|-\|x_{n+1}-p\|^{2}+d_{n}\\
&\leq&\|x_{n}-x_{n+1}\|(\|x_{n}-p\|+\|x_{n+1}-p\|)+d_{n},\\
\end{eqnarray*}
by using (\ref{234}), we have, $\lim_{n\to\infty}\|x_{n}-u_{n}\|=0.$\\
Now we can say that,
$\|y_{n}-x_{n}\|\leq\|y_{n}-u_{n}\|+\|u_{n}-x_{n}\|,$ \\
so,\\
\begin{equation}\label{345}
\lim_{n\to\infty}\|y_{n}-x_{n}\|=0.
\end{equation}

\noindent Next we will show that, $\lim\|Wt_{n}-t_{n}\|\to0$. \\
First we prove that, $\lim_{n\to\infty}\|Bz_{n}-Bp\|=0$.
\begin{eqnarray*}
\|y_{n}-p\|^{2}&=&\|\beta_{n}\gamma
f(z_{n})+(I-\beta_{n}A)t_{n}-p\|^{2}\\
&=&\|\beta_{n}(\gamma f(z_{n})-Ap)+(I-\beta_{n}A)(t_{n}-p)\|^{2}\\
&=&(1-\beta_{n}\bar{\gamma})^{2}\|t_{n}-p\|^{2}+\beta_{n}^{2}\|\gamma
f(z_{n})-Ap\|^{2}\\
&&+2\langle \beta_{n}(\gamma
f(z_{n})-Ap),(I-\beta_{n}A)(t_{n}-p)\rangle\\
&\leq&(1-\beta_{n}\bar{\gamma})^{2}\|z_{n}-\gamma_{n}Bz_{n}-p+\gamma_{n}Bp\|^{2}+\beta_{n}^{2}\|\gamma
f(z_{n})-Ap\|^{2}\\
&&+2\langle \beta_{n}(\gamma
f(z_{n})-Ap),(I-\beta_{n}A)(t_{n}-p)\rangle\\
&=&(1-\beta_{n}\bar{\gamma})^{2}[\|z_{n}-p\|^{2}+\gamma_{n}(\gamma_{n}-2\alpha)\|Bz_{n}-Bp\|^{2}]+c_{n},\\
\end{eqnarray*}
where, $c_{n}=\beta_{n}^{2}\|\gamma f(z_{n})-Ap\|^{2}+2\langle
\beta_{n}(\gamma f(z_{n})-Ap),(I-\beta_{n}A)(t_{n}-p)\rangle.$\\
So,\\
\begin{eqnarray*}
(1-\beta_{n}\bar{\gamma})^{2}b(b-2\alpha)\|Bz_{n}-Bp\|^{2}&\leq&(1-\beta_{n}\bar{\gamma})\|u_{n}-p\|^{2}-\|y_{n}-p\|+c_{n}\\
&\leq&1-\beta_{n}\bar{\gamma})^{2}(\|u_{n}-y_{n}\|+\|y_{n}-p\|)^{2}-\|y_{n}-p\|^{2}+c_{n}\\
&\leq&(1-\beta_{n}\bar{\gamma})^{2}[\|u_{n}-y_{n}\|^{2}+2\|y_{n}-u_{n}\|\|y_{n}-p\|+\|y_{n}-p\|^{2}]\\
&&-\|y_{n}-p\|^{2}+c_{n},\\
\end{eqnarray*}
hence, $\lim_{n\to\infty}\|Bz_{n}-Bp\|=0$. \\
On the other hand,\\
\begin{eqnarray*}
\|t_{n}-p\|^{2}&=&\|P_{C}(z_{n}-\gamma
Bz_{n})-P_{C}(p-\gamma_{n}Bp)\|^{2}\\
&\leq&\langle(z_{n}-\gamma_{n}Bz_{n})-(p-\gamma_{n}Bp),t_{n}-p\rangle\\
&=&\frac{1}{2}\{\|(z_{n}-\gamma_{n}Bz_{n})-(p-\gamma_{n}Bp)\|^{2}+\|t_{n}-p\|^{2}\\
&&-\|(z_{n}-\gamma_{n}Bz_{n})-(p-\gamma_{n}Bp)-(t_{n}-p)\|^{2}\}\\
&\leq&\frac{1}{2}\{\|z_{n}-p\|^{2}+\|t_{n}-p\|^{2}-\|(z_{n}-t_{n})-\gamma_{n}(Bz_{n}-Bp)\|^{2}\}\\
&=&\frac{1}{2}\{\|z_{n}-p\|^{2}+\|t_{n}-p\|^{2}-\|z_{n}-t_{n}\|^{2}\\
&&+2\gamma_{n}\langle
z_{n}-t_{n},Bz_{n}-Bp\rangle-\gamma_{n}^{2}\|Bz_{n}-Bp\|^{2}\},\\
\end{eqnarray*}
which implies,\\
$\|t_{n}-p\|^{2}\leq\|z_{n}-p\|^{2}-\|z_{n}-t_{n}\|^{2}+2\gamma_{n}\langle
z_{n}-t_{n},Bz_{n}-Bp\rangle-\gamma_{n}^{2}\|Bz_{n}-Bp\|^{2}.$\\

\noindent Now by above inequality we want to show that,
$\|z_{n}-t_{n}\|\to0$.
\begin{eqnarray*}
\|y_{n}-p\|^{2}&=&(1-\beta_{n}\bar{\gamma})^{2}\|t_{n}-p\|^{2}+c_{n}\\
&\leq&(1-\beta_{n}\bar{\gamma})^{2}(\|z_{n}-p\|^{2}-\|z_{n}-t_{n}\|^{2}+2\gamma_{n}\langle
z_{n}-t_{n},Bz_{n}-Bp\rangle\\
&&-\gamma_{n}^{2}\|Bz_{n}-Bp\|^{2})+c_{n}.\\
\end{eqnarray*}
Therefore,\\
\begin{eqnarray*}
(1-\beta_{n}\bar{\gamma})^{2}\|z_{n}-t_{n}\|^{2}&\leq&(1-\beta_{n}\bar{\gamma})^{2}\|z_{n}-p\|^{2}-\|y_{n}-p\|^{2}\\
&&+2\gamma_{n}\|z_{n}-t_{n}\|\|Bz_{n}-Bp\|-\gamma_{n}^{2}\|Bz_{n}-Bp\|^{2}+c_{n}\\
&\leq&(1-\beta_{n}\bar{\gamma})^{2}\|u_{n}-p\|^{2}-\|y_{n}-p\|^{2}+2\gamma_{n}\|z_{n}-t_{n}\|\|Bz_{n}-Bp\|\\
&&-\gamma_{n}^{2}\|Bz_{n}-Bp\|^{2}+c_{n}\\
&\leq&(1-\beta_{n}\bar{\gamma})^{2}[\|u_{n}-y_{n}\|^{2}+2\|y_{n}-u_{n}\|\|y_{n}-p\|+\|y_{n}-p\|^{2}]-\|y_{n}-p\|^{2}\\
&&+2\gamma_{n}\|z_{n}-t_{n}\|\|Bz_{n}-Bp\|-\gamma_{n}^{2}\|Bz_{n}-Bp\|^{2}+c_{n},\\
\end{eqnarray*}
so, $\lim_{n\to\infty}\|z_{n}-t_{n}\|=0$.\\
\begin{equation}\label{444}
\end{equation}
Also,\\
\begin{eqnarray*}
\|W_{n}y_{n}-y_{n}\|&\leq&\|y_{n}-z_{n}\|+\|z_{n}-W_{n}y_{n}\|\\
&\leq&\|y_{n}-z_{n}\|+\|z_{n}-W_{n}u_{n}\|+\|W_{n}u_{n}-W_{n}y_{n}\|\\
&\leq&\|y_{n}-z_{n}\|+\lambda_{n}\|u_{n}-W_{n}u_{n}\|+\|u_{n}-y_{n}\|\\
&\leq&\|y_{n}-z_{n}\|+\lambda_{n}\|u_{n}-y_{n}\|+\lambda_{n}\|y_{n}-W_{n}u_{n}\|+\|u_{n}-y_{n}\|\\
&\leq&\|y_{n}-z_{n}\|+(1+\lambda_{n})\|u_{n}-y_{n}\|+\lambda_{n}\|y_{n}-W_{n}y_{n}\|\\
&&+\lambda_{n}\|W_{n}y_{n}-W_{n}u_{n}\|.\\
\end{eqnarray*}
So,\\
$(1-\lambda_{n})\|W_{n}y_{n}-y_{n}\|\leq\|y_{n}-z_{n}\|+(2+\lambda_{n})\|u_{n}-y_{n}\|,$ by (\ref{444}) and (\ref{333}) we have,\\
\begin{equation}\label{mh}
\|W_{n}y_{n}-y_{n}\|\to0.\\
\end{equation}
By lemma (\ref{2568}) and (\ref{mh}), we obtain,
\begin{eqnarray*}
\|Wy_{n}-y_{n}\|&\leq&\|Wy_{n}-W_{n}y_{n}\|+\|W_{n}y_{n}-y_{n}\|\\
&\leq&\sup_{y\in \{y_{n}\}}\|Wy-W_{n}y\|+\|W_{n}y_{n}-y_{n}\|\to0, \
\ as \ \ \ n\to\infty.\\
\end{eqnarray*}
\noindent Finally we have,\\
$\|Wt_{n}-t_{n}\|\leq\|Wt_{n}-Wy_{n}\|+\|Wy_{n}-y_{n}\|+\|y_{n}-t_{n}\|.$
So, $\|Wt_{n}-t_{n}\|\to0$.
\begin{equation}\label{2434}
\end{equation}

\noindent Observe that, $P_{F}(\gamma f+(I-A))$ is a contraction, in fact for $x, y\in H$ we have,\\
\begin{eqnarray*}
\|P_{F}(\gamma f+(I-A))(x)-P_{F}(\gamma f+(I-A))(y)\|&\leq&\|(\gamma f+(I-A))(x)-(\gamma f+(I-A))(y)\|\\
&\leq&\gamma\|f(x)-f(y)\|+\|I-A\|\|x-y\|\\
&\leq&\gamma\alpha\|x-y\|+(1-\bar{\gamma})\|x-y\|\\
&=&[1-(\bar{\gamma}-\gamma\alpha)]\|x-y\|.\\
\end{eqnarray*}
Banach's Contraction Mapping Principle guarantees that $P_{F}(\gamma f+(I-A))$ has a unique fixed point, say $q\in H$. That is, $P_{F}(\gamma f+(I-A))(q)=q$.\\
Now we choose a subsequence $\{t_{n_{i}}\}$ of $\{t_{n}\}$ such that,\\
$$\limsup_{n\to\infty}\langle\gamma f(q)-Aq,Wt_{n}-q\rangle=\lim_{i\to\infty}\langle\gamma f(q)-Aq,Wt_{n_{i}}-q\rangle.$$
Since $\{t_{n_{i}}\}$ is bounded, there exists a subsequence
$\{t_{n_{i_{j}}}\}$ of $\{t_{n_{i}}\}$ which converges weakly to
$z\in C$. Without loss of generality,
we can assume that $t_{n_{i}}\rightharpoonup z$. From $\|Wt_{n_{i}}-t_{n_{i}}\|\to0$, we obtain $Wt_{n_{i}}\rightharpoonup z$. Therefore, we have \\
\begin{eqnarray*}
\limsup_{n\to\infty}\langle\gamma f(q)-Aq,Wt_{n}-q\rangle&=&\lim_{i\to\infty}\langle\gamma f(q)-Aq,Wt_{n_{i}}-q\rangle\\
&=&\langle\gamma f(q)-Aq,z-q\rangle.\\
\end{eqnarray*}
\begin{equation}\label{sha}
\end{equation}

\noindent We want to prove that $z\in
F(W)=\bigcap_{i=1}^{\infty}F(T_{i})$.
Suppose $z\notin F(W)$, that is $Wz\neq z$. Since $t_{n_{i}}\rightharpoonup z$, by Opial's condition and (\ref{2434}) we have,\\
\begin{eqnarray*}
\liminf_{i\to\infty}\|t_{n_{i}}-z\|&<&\liminf_{i\to\infty}\|t_{n_{i}}-Wz\|\\
&\leq&\liminf_{i\to\infty}\{\|t_{n_{i}}-Wt_{n_{i}}\|+\|Wt_{n_{i}}-Wz\|\}\\
&\leq&\liminf_{i\to\infty}\{\|t_{n_{i}}-Wt_{n_{i}}\|+\|t_{n_{i}}-z\|\}\\
&=&\liminf_{i\to\infty}\|t_{n_{i}}-z\|.\\
\end{eqnarray*}
That is a contradiction, hence $z\in F(W)=\bigcap_{i=1}^{\infty}F(T_{i})$.\\
Next we prove that $z\in EP(F)$. Since $u_{n}=T_{r_{n}}$, we have,\\
$$F(u_{n},y)+\frac{1}{r_{n}}\langle y-u_{n},u_{n}-x_{n}\rangle\geq0, \ \forall y\in C.$$
From (A2), we have, $$\frac{1}{r_{n}}\langle
y-u_{n},u_{n}-x_{n}\rangle \geq F(y,u_{n})$$
and hence,\\
$$\langle y-u_{n_{i}},\frac{u_{n_{i}}-x_{n_{i}}}{r_{n_{i}}}\rangle\geq F(y,u_{n_{i}}).$$
From $\|u_{n}-x_{n}\|\to0$, $\|x_{n}-Wt_{n}\|\to0$ and
$\|Wt_{n}-t_{n}\|\to0$ we get $u_{n_{i}}\rightharpoonup z$.Since
$\frac{u_{n_{i}}-x_{n_{i}}}{r_{n_{i}}}\to0$
by (A4) we have, $F(y,z)\leq0, \ \forall y\in C$.\\
For $t$ with $0<t\leq 1$ and $y\in C$, let $y_{t}=ty+(1-t)z$. Since
$y\in C$ and $z\in C$ we have $y_{t}\in C$ and hence
$F(y_{t},z)\leq0$. So from (A1) and (A4) we have
$$0=F(y_{t},y_{t})\leq tF(y_{t},y)+(1-t)F(y_{t},z)\leq tF(y_{t},y)$$
and hence $0\leq F(y_{t},y)$. From (A3) we have $0\leq F(z,y)$ for all $y\in C$ so, $z\in EP(F)$.\\

\noindent Finally, similar as the proof of [9, Theorem 2.1, pp. 678], we can prove that $z\in VI(B,C)$. Therefore $$z\in F(W)\cap VI(B,C)\cap EP(F).$$\\
Moreover by (\ref{sha}), we have,\\
$$\limsup_{n\to\infty}\langle \gamma f(q)-Aq,Wt_{n}-q\rangle=\langle \gamma
f(q)-Aq,z-q\rangle\leq0.$$ Also,
\begin{eqnarray*}
\limsup_{n\to\infty}\langle \gamma f(q)-Aq,t_{n}-q\rangle&=&\limsup_{n\to\infty}\langle \gamma f(q)-Aq,(t_{n}-Wt_{n})+(Wt_{n}-q)\rangle\\
&\leq&\limsup_{n\to\infty}\langle \gamma f(q)-Aq,Wt_{n}-q\rangle\\
&\leq&0.\\
\end{eqnarray*}
\begin{equation}\label{2222}
\end{equation}
Hence we have,\\
\begin{eqnarray*}
\limsup_{n\to\infty}\langle \gamma f(q)-Aq,x_{n}-q\rangle&=&\langle \gamma f(q)-Aq,x_{n}-y_{n}\rangle\\
&&+\langle \gamma f(q)-Aq,y_{n}-t_{n}\rangle+\langle \gamma f(q)-Aq,t_{n}-q\rangle,\\
\end{eqnarray*}
by (\ref{345}) and (\ref{2222}) we have,  $\limsup_{n\to\infty}\langle \gamma f(q)-Aq,x_{n}-q\rangle\leq0$.\\
\begin{equation}\label{rq}
\end{equation}

\noindent On the other hand we have, $\langle \gamma
f(q)-Aq,y_{n}-q\rangle=\langle \gamma
f(q)-Aq,y_{n}-x_{n}\rangle+\langle\gamma f(q)-Aq,x_{n}-q\rangle,$\\
by (\ref{345}) and (\ref{rq}) we have,
\begin{equation}\label{567} \limsup_{n\to\infty}\langle \gamma
f(q)-Aq,y_{n}-q\rangle\leq0.\\
\end{equation}
So,\\
\begin{eqnarray*}
\|y_{n}-q\|^{2}&=&\|\beta_{n}\gamma
f(z_{n})+(I-\beta_{n}A)t_{n}-q\|^{2}\\
&=&\|(I-\beta_{n}A)(t_{n}-q)+\beta_{n}(\gamma f(z_{n})-Aq)\|^{2}\\
&\leq&\|(I-\beta_{n}A)(t_{n}-q)\|^{2}+2\beta_{n}\langle \gamma
f(z_{n})-Aq,y_{n}-q\rangle\\
&=&\|(I-\beta_{n}A)(t_{n}-q)\|^{2}+2\beta_{n}\langle \gamma
f(z_{n})-\gamma f(q),y_{n}-q\rangle\\
&&+2\beta_{n}\langle \gamma f(q)-Aq,y_{n}-q\rangle\\
&\leq&(1-\beta_{n}\bar{\gamma})^{2}\|t_{n}-q\|^{2}+2\beta_{n}\gamma\alpha\|z_{n}-q\|\|y_{n}-q\|\\
&&+2\beta_{n}\langle \gamma f(q)-Aq,y_{n}-q\rangle\\
&\leq&(1-\beta_{n}\bar{\gamma})^{2}\|x_{n}-q\|^{2}+2\beta_{n}\gamma\alpha\|x_{n}-q\|\|y_{n}-q\|\\
&&+2\beta_{n}\langle \gamma f(q)-Aq,y_{n}-q\rangle\\
&\leq&(1-\beta_{n}\bar{\gamma})^{2}\|x_{n}-q\|^{2}+\beta_{n}\gamma\alpha(\|x_{n}-q\|^{2}+\|y_{n}-q\|^{2})\\
&&+2\beta_{n}\langle\gamma f(q)-Aq,y_{n}-q\rangle,\\
\end{eqnarray*}
so,\\
\begin{eqnarray*}
\|y_{n}-q\|^{2}&\leq&\frac{(1-\beta_{n}\bar{\gamma})^{2}+\beta_{n}\gamma\alpha}{1-\beta_{n}\gamma\alpha}\|x_{n}-q\|^{2}\\
&&+\frac{2\beta_{n}}{1-\beta_{n}\gamma\alpha}\langle \gamma f(q)-Aq,y_{n}-q\rangle\\
&\leq&[1-\frac{2\beta_{n}(\bar{\gamma}-\gamma\alpha)}{1-\beta_{n}\gamma\alpha}]\|x_{n}-q\|^{2}\\
&&+\frac{2\beta_{n}(\bar{\gamma}-\gamma\alpha)}{1-\beta_{n}\gamma\alpha}[\frac{1}{\bar{\gamma}-\gamma\alpha}\langle
\gamma
f(q)-Aq,y_{n}-q\rangle+\frac{\beta_{n}\bar{\gamma}^{2}}{2(\bar{\gamma}-\gamma\alpha)}M_{4}],\\
\end{eqnarray*}
where $M_{4}$ is an appropriate constant such that
$M_{4}\geq\sup_{n\geq1}\{\|x_{n}-q\|\}.$

\noindent Finally,\
\begin{eqnarray*}
\|x_{n+1}-q\|^{2}&\leq&\alpha_{n}\|x_{n}-q\|^{2}+(1-\alpha_{n})\|y_{n}-q\|^{2}\\
&\leq&[1-(1-\alpha_{n})\frac{2\beta_{n}(\bar{\gamma}-\gamma\alpha)}{1-\beta_{n}\gamma\alpha}\|x_{n}-q\|^{2}\\
&&+(1-\alpha_{n})\frac{2\beta_{n}(\bar{\gamma}-\gamma\alpha)}{1-\beta_{n}\gamma\alpha}[\frac{1}{\bar{\gamma}-\gamma\alpha}\langle
\gamma
f(q)-Aq,y_{n}-q\rangle+\frac{\beta_{n}\bar{\gamma}}{2(\bar{\gamma}-\gamma\alpha)}M_{4}],\\
\end{eqnarray*}
we put
$h_{n}=(1-\alpha_{n})\frac{2\beta_{n}(\bar{\gamma}-\gamma\alpha)}{1-\beta_{n}\gamma\alpha}$
and $r_{n}=\frac{1}{\bar{\gamma}-\gamma\alpha}\langle \gamma
f(q)-Aq,y_{n}-q\rangle+\frac{\beta_{n}\bar{\gamma}}{2(\bar{\gamma}-\gamma\alpha)}M_{4},$\\
so,\\
\begin{equation}\label{456}
\|x_{n+1}-q\|^{2}\leq(1-h_{n})\|x_{n}-q\|^{2}+t_{n},
\end{equation}
where $t_{n}=h_{n}r_{n}$, from condition (C1) we have
$\lim_{n\to\infty}h_{n}=0$ and by(\ref{567}),
$\limsup_{n\to\infty}t_{n}\leq0$. \\
Finally lemma (\ref{1}), implies $x_{n}\to0$ as $n\to\infty$.  \   \
 \   \    \    \  \ \
\   \   \   \  \  \    \ \  \ \    \ \    \  \    \   \    \   \ \   \ \  \   $\square$\\

\noindent If we put $T_{n}x=x$ for all $n=1,2,\ldots,k$ and for all
$x\in C$ in (\ref{12}) we have $W_{n}x=x$, then we have,\\
\begin{cor}
Let $H$ be a real Hilbert space. Let $F$ be a bifunction from
$C\times C \longrightarrow \mathbb{R}$ satisfying (A1)-(A4) and let
$B$ be an $\alpha$-inverse-strongly monotone mapping of $C$ into
$H$, $F:=VI(B,C)\bigcap EP(F)\neq \emptyset$. Suppose $A$ is a
strongly positive linear bounded self-adjoint operator with the
coefficient $\bar{\gamma}\geq0$ and
$0<\gamma\leq\frac{\bar{\gamma}}{\alpha}$. Let $\{\alpha_{n}\}$,
$\{\beta_{n}\}$, $\{\lambda_{n}\}$ and $\{\delta_{n}\}$ be sequences
in $[0,1]$ satisfying the following conditions:\\
(C1) $\sum_{n=0}^{\infty}\beta_{n}=\infty$,
$\lim_{n\to\infty}\beta_{n}=0$;\\
(C2) $\sum_{i=1}^{\infty}|\lambda_{n}-\lambda_{n+1}|<\infty$;\\
(C3)
$\sum_{i=1}^{\infty}|\alpha_{n}-\alpha_{n+1}|<\infty;$\\
(C4) $\exists \lambda\in[0,1]; \lambda_{n}<\lambda, \ \forall
n\geq1;$\\
(C5) $\sum_{n=1}^{\infty}|\gamma_{n}-\gamma_{n+1}|<\infty,
\gamma_{n}\in[a,b], a,b\in(0,2\alpha);$\\
(C6) $\liminf_{n\to\infty}r_{n}>0$ and $\sum_{i=1}^{\infty}|r_{n+1}-r_{n}|<\infty$.\\
Then the sequences $\{x_{n}\}$ defined by
\begin{eqnarray*}
&&F(u_{n},y)+\frac{1}{r_{n}}\langle y-u_{n},u_{n}-x_{n}\rangle\geq0,
\
\forall y\in C, \\
&&y_{n}=\beta_{n}\gamma
f(u_{n})+(I-\beta_{n}A)P_{C}(u_{n}-\gamma_{n}Bu_{n}),  \\
&&x_{n+1}=\alpha_{n}x_{n}+(1-\alpha_{n})y_{n},  \\
\end{eqnarray*}
converges strongly to $q\in F$, where $q=P_{F}(\gamma f+(I-A))(q)$
which solves the following variational inequality: $$\langle \gamma
f(q)-Aq,p-q\rangle\leq0, \ \forall p\in F.$$

\end{cor}

\noindent If we put $F(x,y)=0$ for all $x, y \in C$ and $r_{n}=1$
for all $n\in\mathbb{N}$ in theorem (\ref{maintheorem}), then we
have,
\begin{cor}
Let $H$ be a real Hilbert space. Let $\{T_{i}:C\longrightarrow C\}$
be a infinite family of nonexpansive mappings with
$F:=\bigcap_{i=1}^{\infty}F(T_{i})\bigcap VI(B,C)\neq \emptyset$.
Suppose $A$ is a strongly positive linear bounded self-adjoint
operator with the coefficient $\bar{\gamma}\geq0$ and
$0<\gamma\leq\frac{\bar{\gamma}}{\alpha}$. Let $\{\alpha_{n}\}$,
$\{\beta_{n}\}$, $\{\lambda_{n}\}$ and $\{\delta_{n}\}$ be sequences
in $[0,1]$ satisfying the following conditions:\\
(C1) $\sum_{n=0}^{\infty}\beta_{n}=\infty$,
$\lim_{n\to\infty}\beta_{n}=0$;\\
(C2) $\sum_{i=1}^{\infty}|\lambda_{n}-\lambda_{n+1}|<\infty$;\\
(C3)
$\sum_{i=1}^{\infty}|\alpha_{n}-\alpha_{n+1}|<\infty;$\\
(C4) $\exists \lambda\in[0,1]; \lambda_{n}<\lambda, \ \forall
n\geq1;$\\
(C5) $\sum_{n=1}^{\infty}|\gamma_{n}-\gamma_{n+1}|<\infty,
\gamma_{n}\in[a,b], a,b\in(0,2\alpha);$\\
Then the sequences $\{x_{n}\}$ defined by
\begin{eqnarray*}
&&z_{n}=\lambda_{n}u_{n}+(1-\lambda_{n})W_{n}u_{n},\\
&&y_{n}=\beta_{n}\gamma
f(z_{n})+(I-\beta_{n}A)P_{C}(z_{n}-\gamma_{n}Bz_{n}), \\
&&x_{n+1}=\alpha_{n}x_{n}+(1-\alpha_{n})y_{n}, \\
\end{eqnarray*}
 converges strongly to $q\in
F$, where $q=P_{F}(\gamma f+(I-A))(q)$ which solves the following
variational inequality: $$\langle \gamma f(q)-Aq,p-q\rangle\leq0, \
\forall p\in F.$$

\end{cor}

\end{document}